\magnification=\magstep1
\input amstex
\documentstyle{amsppt}

\define\defeq{\overset{\text{def}}\to=}
\define\ab{\operatorname{ab}}
\define\pr{\operatorname{pr}}
\define\Gal{\operatorname{Gal}}
\define\Ker{\operatorname{Ker}}
\define\genus{\operatorname{genus}}
\define\et{\operatorname{et}}

\define\id{\operatorname{id}}
\define\Pic{\operatorname{Pic}}
\def \isom {\overset \sim \to \rightarrow}
\define\Sec{\operatorname{Sec}}
\define\Spec{\operatorname{Spec}}

\def\nor{\operatorname{nor}}

\NoRunningHeads
\NoBlackBoxes
\topmatter

\title
On the existence of non-geometric sections of arithmetic fundamental groups
\endtitle

\author
Mohamed Sa\"\i di
\endauthor
\abstract
We show the existence of group-theoretic sections of the {\it "\'etale-by-geometrically abelian" quotient} 
of the arithmetic fundamental group of hyperbolic curves over $p$-adic local fields {\it relative to a proper and flat model} which are {\it non-geometric}, i.e.,
which do not arise from rational points.
\endabstract

\endtopmatter

\document

\subhead
\S 0. Introduction/Statement of the Main Result
\endsubhead
This note is motivated by the {\it $p$-adic analog of the Grothendieck anabelian section conjecture}, which asks
if {\it group-theoretic sections of arithmetic fundamental groups of hyperbolic curves over $p$-adic local fields
all arise from rational points}. (See [Sa\"\i di] for more details on the statement of this conjecture).  
For the time being it is not known if this conjecture holds or not (cf. loc. cit. for a conditional proof of this conjecture). 

Recently, examples were found of group-theoretic sections of certain {\it quotients of 
arithmetic fundamental groups} of curves over $p$-adic local fields which are {\it non-geometric}, i.e.,
which {\it do not arise from rational points}. Hoshi has provided examples of sections of the {\it geometrically pro-$p$} quotient
of arithmetic fundamental groups of curves over $p$-adic local fields which are non-geometric (cf. [Hoshi]).
(Actually, Hoshi's example arises from group-theoretic sections of geometrically pro-$p$ fundamental groups
of hyperbolic curves over number fields (cf. loc. cit.)). In [Sa\"\i di1]  we provided examples of group-theoretic sections of 
{\it geometrically prime-to-$p$} fundamental groups
of hyperbolic curves over $p$-adic local fields which are non-geometric (cf. loc. cit. $\S3$). 
{\it The existence of these examples is crucial for our understanding 
of the $p$-adic section conjecture}. Indeed, if the $p$-adic version of the section conjecture
holds true then it may possibly hold true even for smaller quotients of the arithmetic fundamental group,
and one would like to know these quotients in this case. On the other hand, more elaborate examples
of non-geometric sections as above may lead to a counterexample for the $p$-adic version of the section conjecture.

In this note {\it we provide examples of sections of certain quotients of arithmetic fundamental groups of 
curves over $p$-adic local fields which are non-geometric}. 
(Our examples, including the quotients of arithmetic fundamental groups involved,
are quite different from those considered in  [Hoshi]). The quotients of arithmetic fundamental groups 
we consider are roughly speaking the "{\it \'etale-by-geometrically abelian}" quotient of the arithmetic fundamental group of
a curve over a $p$-adic local field {\it relative to a proper and flat model} over the ring of integers of the base field 
(cf. the discussion before Theorem A for a precise definition of this quotient).

Next, we fix some notations and state our main results.
Let 
$$1\to H'\to H @>{\pr}>> G \to 1$$ 
be an exact sequence of profinite groups. We will refer to
a continuous homomorphism $s:G\to H$ such that $\pr \circ s=\id_{G}$ as a (group-theoretic)
{\bf section} of the above sequence, or simply a section of the projection $\pr : H\twoheadrightarrow G$.

In this paper $p\ge 2$ is a {\it prime integer}. 
Let $k$ be a $p$-{\bf adic local field}; i.e., $k/\Bbb Q_p$ is a finite extension,
with ring of integers $\Cal O_k$, and residue field $F$. 
Thus, $F$ is a {\it finite field} of characteristic $p$. 
Let $X\to \Spec k$ 
be a proper, smooth, and geometrically connected {\bf hyperbolic} (i.e., $\genus (X)\ge 2$) {\bf curve} over $k$,
and $\Cal X\to \Spec \Cal O_k$
a {\bf proper}, and {\bf flat model} of $X$ over $\Cal O_k$. We have a commutative diagram with cartesian squares
$$
\CD
X  @>>>   \Spec k\\
@VVV      @VVV\\
\Cal X  @>>> \Spec \Cal O_k\\
@AAA          @AAA\\
\Cal X_s  @>>> \Spec F\\
\endCD
$$
where the right vertical maps are the natural ones, and $\Cal X_s\defeq \Cal X\times _{\Cal O_k}F$ is {\it the special fibre} of $\Cal X$.

Let $\eta$ be a geometric point of $X$ above the generic point of $X$.
Then $\eta$ determines naturally an algebraic closure $\overline k$
of $k$, and a geometric point $\bar {\eta}$ of $\overline {X} \defeq X\times _k \overline k$.
There exists a canonical exact sequence of profinite groups (cf. [Grothendieck], Expos\'e IX, Th\'eor\`eme 6.1)
$$1\to \pi_1(\overline {X},\bar \eta)\to \pi_1(X, \eta) @>>> G_k\to 1.$$
Here, $\pi_1(X, \eta)$ denotes the {\it arithmetic \'etale fundamental group} of $X$ with base
point $\eta$, $\pi_1(\overline {X},\bar \eta)$ the \'etale fundamental group of $\overline {X}\defeq 
X\times _k \overline k$ with base
point $\bar \eta$, and $G_k\defeq \Gal (\overline k/k)$ the absolute Galois group of $k$.

Let $\xi$ be a geometric point of $\Cal X_s$ above the generic point of some {\it irreducible component} $X_{i_0}$ of $\Cal X_s$.
Then $\xi$ determines naturally an algebraic closure $\overline F$
of $F$, and a geometric point $\bar {\xi}$ of $\overline {\Cal X_s} \defeq \Cal X_s\times _F \overline F$.
There exists a canonical exact sequence of profinite groups (cf. loc. cit.)
$$1\to \pi_1(\overline {\Cal X_s},\bar \xi)\to \pi_1(\Cal X_s, \xi) @>>> G_F\to 1.$$
Here, $\pi_1(\Cal X_s, \xi)$ denotes the {\it arithmetic \'etale fundamental group} of $\Cal X_s$ with base
point $\xi$, $\pi_1(\overline {\Cal X_s},\bar \xi)$ the \'etale fundamental group of $\overline {\Cal X_s}\defeq 
\Cal X_s\times _F \overline F$ with base
point $\bar \xi$, and $G_F\defeq \Gal (\overline F/F)$ the absolute Galois group of $F$.

In what follows we {\it assume that the $\gcd$ of the multiplicities of the irreducible components of $\Cal X_s$ equals $1$}.
Then, after a suitable choice of the base points $\xi$ and $\eta$, there exists a natural commutative 
diagram

$$
\CD
@.  1@.   1@.   1@.\\
@. @VVV   @VVV  @VVV \\
 1 @>>>  I_{\overline {\Cal X}}   @>>>   I_{\Cal X}  @>>>   I_k  @>>> 1\\
@.  @VVV     @VVV    @VVV  \\
1 @>>>    \pi_1(\overline {X},\bar \eta)   @>>>  \pi_1(X, \eta)      @>>> G_k @>>> 1\\
  @.        @VVV           @VVV              @VVV \\
1 @>>>   \pi_1(\overline {\Cal X_s},\bar \xi)    @>>>   \pi_1(\Cal X_s, \xi)       @>>> G_F   @>>> 1\\
@. @VVV   @VVV  @VVV \\
@.  1@.    1@.   1@.\\
\endCD
\tag 0.1
$$
where the vertical and horizontal sequences are exact, $I_k$ (resp. $I_{\Cal X}$, and $I_{\overline {\Cal X}}$)
are defined so that the right (resp. the middle, and left) vertical sequence is exact 
($I_k$ is the inertia subgroup of $G_k$);
as follows easily from the {\it specialisation theory} for fundamental groups
of Grothendieck (cf. [Grothendieck], Expos\'e X, $\S2$, and [Raynaud], Proposition 6.3.5, for the surjectivity of the left lower vertical map under our assumption
on the $\gcd$ of the multiplicities of the components of $\Cal X_s$, and the fact that $F$ is perfect). 
Also,  consider the following commutative diagram

$$
\CD
1 @>>>    \pi_1(\overline {X},\bar \eta)   @>>>  \pi_1(X, \eta)      @>>> G_k @>>> 1\\
  @.        @VVV           @VVV              @VVV \\
1 @>>>   \pi_1(\overline {\Cal X_s},\bar \xi)    @>>>   \widehat \pi_1(\Cal X_s, \xi)       @>>> G_k  @>>> 1\\
@. @VVV   @VVV  @VVV \\
1 @>>>   \pi_1(\overline {\Cal X_s},\bar \xi)    @>>>   \pi_1(\Cal X_s, \xi)       @>>> G_F   @>>> 1\\
@. @VVV   @VVV  @VVV \\
@.  1@.    1@.   1@.\\
\endCD
\tag 0.2
$$
where the horizontal sequences are exact, the upper and lower horizontal sequences are the middle and lower, respectively, horizontal 
sequences in diagram (0.1),
the lower right square is {\it cartesian}, and all vertical maps are surjective.
Thus, $\widehat \pi_1(\Cal X_s, \xi)$ is the {\it pullback} of $\pi_1(\Cal X_s, \xi)$ via the natural projection
$G_k\twoheadrightarrow G_F$. We shall refer to the quotient 
$$\pi_1(X,\eta)\twoheadrightarrow  \widehat \pi_1(\Cal X_s, \xi)$$  
of $\pi_1(X,\eta)$ as the {\bf \'etale} quotient of $\pi_1(X,\eta)$ {\bf relative to the model $\Cal X$}. 
We have an exact sequence (cf. diagram (0.1), and diagram (0.2))
$$1\to I_{\overline {\Cal X}} \to \pi_1(X,\eta)\to  \widehat \pi_1(\Cal X_s, \xi)\to 1.$$

Write
$I_{\overline {\Cal X}}^{\ab}$ for the maximal {\it abelian} quotient of $I_{\overline {\Cal X}}$, which is a {\it characteristic} quotient. 
Consider the following {\it pushout} diagrams

$$
\CD
1 @>>>  I_{\overline {\Cal X}} @>>>  \pi_1(\overline {X},\bar \eta) @>>>   \pi_1(\overline {\Cal X_s},\bar \xi)  @>>> 1 \\
@.      @VVV     @VVV   @VVV  \\
1 @>>> I_{\overline {\Cal X}}^{\ab}    @>>>  \Delta^{\et, \ab}_{\overline {\Cal X}} @>>>    \pi_1(\overline {\Cal X_s},\bar \xi)  @>>> 1 \\
\endCD
\tag 0.3.1
$$
which defines a natural quotient  $\Delta ^{\et,\ab}_{\overline {\Cal X}}$ of $\pi_1(\overline {X},\bar \eta)$, 
and 

$$
\CD
1 @>>>  I_{\overline {\Cal X}} @>>>  I_{\Cal X} @>>>  I_k @>>> 1 \\
@.      @VVV     @VVV   @VVV  \\
1 @>>> I_{\overline {\Cal X}}^{\ab}    @>>>  I_{\Cal X}^{(\ab)} @>>>    I_k  @>>> 1 \\
\endCD
\tag 0.3.2
$$
which defines a natural quotient $I_{\Cal X}^{(\ab)}$ of $I_{\Cal X}$. Let 
$$R_{\Cal X}\defeq \Ker (I_{\Cal X}\to I_{\Cal X}^{(\ab)})=\Ker (I_{\overline {\Cal X}}\to I_{\overline {\Cal X}}^{\ab})=
\Ker (\pi_1(\overline X,\bar \eta)\to \Delta^{\et, \ab}_{\Cal X}).$$ 
Note that $R_{\Cal X}$ is a {\it normal subgroup} of $\pi_1(X, \eta)$. Indeed, this follows from the fact that $I_{\overline {\Cal X}}$ 
is the kernel of the natural surjective homomorphism $\pi_1(X,\eta)\twoheadrightarrow \hat \pi_1(\Cal X_s,\xi)$, 
and $I_{\overline {\Cal X}}^{\ab}$ is a characteristic quotient of $I_{\overline {\Cal X}}$.
Write 
$$\Pi^{(\et,ab)}_{\Cal X}\defeq \pi_1(X,\eta)/R_{\Cal X}.$$
Note that $I_{\Cal X}^{(\ab)}$ is the image of $I_{\Cal X}$ in $\Pi^{(\et,ab)}_{\Cal X}$ (cf. diagram (0.3.2)).

We have a natural commutative diagram of exact sequences
$$
\CD
1 @>>> \pi_1(\overline {X},\bar \eta) @>>>   \pi_1(X,\eta) @>>>  G_k  @>>> 1 \\
@.      @VVV     @VVV   @VVV  \\
1 @>>> \Delta ^{\et,\ab}_{\overline {\Cal X}}    @>>>  {\Pi}^{(\et, \ab)}_{\Cal X} @>>>    G_k @>>> 1 \\
\endCD
\tag 0.4
$$

We will refer to $\Pi^{(\et,\ab)}_{\Cal X}$ as the {\bf \'etale-by-geometrically abelian} quotient of $\pi_1(X,\eta)$ {\bf relative to the model $\Cal X$}.
Note that a rational point $x\in X(k)$ gives rise naturally (by the fuctoriality of $\pi_1$) to (a conjugacy class of) a  group-theoretic section 
$$s_x:G_k\to \pi_1(X,\eta)$$
of the upper horizontal exact sequence in diagram (0.4), hence gives also rise to a section 
$$s_x:G_k\to  {\Pi}^{(\et, \ab)}_{\Cal X}$$
of the natural projection ${\Pi}^{(\et, \ab)}_{\Cal X}\twoheadrightarrow G_k$
(cf. diagram (0.4)). We will refer to such a section as a {\bf geometric} section of $\Pi^{(\et,\ab)}_{\Cal X}$. 
A {\bf non-geometric} section of the projection $\Pi^{(\et,\ab)}_{\Cal X}\twoheadrightarrow G_k$ is a group-theoretic section 
which doesn't arise from a $k$-rational point $x$ as above.
In this paper we investigate sections of ${\Pi}^{(\et, \ab)}_{\Cal X}\twoheadrightarrow G_k$ 
for a given model $\Cal X$ of $X$. Our main results provide {\it examples of group-theoretic sections of 
the natural projection ${\Pi}^{(\et, \ab)}_{\Cal X}\twoheadrightarrow G_k$ which are non-geometric} (cf. Theorem B, and Theorem C).

\definition {Definition 0.1} We use the same notations as above. Let $\Cal X\to \Spec \Cal O_k$ be a 
{\it proper} and {\it flat} model of $X$ over $\Cal O_k$.
We say that {\it the model $\Cal X$ satisfies the condition $\bold {(\star)}$} if the following holds.

\noindent
(a)\ The scheme $\Cal X$ is {\bf regular}. 

\noindent
(b)\ There exists an {\bf irreducible component} $X_{i_0}$ of the special fibre
$\Cal X_s\defeq \Cal X\times _{\Cal O_k}F$ of $\Cal X$ with the following properties.

\noindent
(i) $X_{i_0}$ is {\bf reduced}. In particular, $X_{i_0}$ is {\it geometrically reduced} since $F$ is perfect,
and {\it the $\gcd$ of the multiplicities of the irreducible components of $\Cal X_s$ equals $1$}.

\noindent
(ii)\  $X_{i_0}$ is {\bf geometrically irreducible}. Thus, $X_{{i_0}}$
is {\it geometrically integral} since $F$ is perfect.

\noindent
(iii)\  $X_{i_0}$ is {\bf geometrically unibranch}, i.e., for each finite extension $F'/F$ the morphism
of normalisation $X_{i_0,F'}^{\nor}\to X_{i_0,F'}\defeq X_{i_0}\times _FF'$ is a {\it homeomorphism}.
\enddefinition

\noindent
(iv)\ The {\bf normalisation} $X_{i_0}^{\nor}$ of $X_{i_0}$ is
{\bf hyperbolic}.

\definition{Remark 0.2} (i)\ If the curve $X$ has {\it good reduction} over $\Cal O_k$, i.e., if $X\to \Spec k$ extends to
a {\it proper and smooth relative curve} $\Cal X\to \Spec \Cal O_k$ over $\Cal O_k$, then the smooth model $\Cal X$ of $X$ 
satisfies the condition $(\star)$.

\noindent
(ii)\ If the curve $X$ satisfies the condition $\bold {(\star)}$, then the {\it index} of $X$ equals $1$ (cf. [Liu], exercise 9.1.9, and the reference therein).
\enddefinition

Our main results in this paper are the following.

\proclaim {Theorem A} Let $X$ be a proper, smooth, and geometrically connected {\bf hyperbolic curve over the $p$-adic local field $k$}.
Let $\Cal X\to \Spec \Cal O_k$ be a proper and flat model of $X$ over $\Cal O_k$.
{\bf Assume that the model $\Cal X$ satisfies the condition} $\bold {(\star)}$ (cf. Definition 0.1). 
Then the corresponding exact sequence
$$1 \to \Delta ^{\et,\ab}_{\overline {\Cal X}}   \to  {\Pi}^{(\et, \ab)}_{\Cal X} \to    G_k \to 1,$$
where $\Pi^{(\et,\ab)}_{\Cal X}$ is the {\bf \'etale-by-geometrically abelian} quotient of the arithmetic fundamental group of $X$ {\bf relative to the model $\Cal X$}
(cf. diagram (0.4)) is a {\bf split} exact sequence of profinite groups.
\endproclaim

As a corollary of Theorem A we obtain the following.

\proclaim {Theorem B} There {\bf exists} a $p$-adic local field $k$, and a smooth, projective, geometrically connected, and {\bf hyperbolic curve} $X$ over
$k$, such that the following holds. There exists a {\bf proper} and {\bf smooth} model $\Cal X\to \Spec \Cal O_k$ of $X$ (i.e., the curve $X$ has 
{\bf good reduction} over $\Cal O_k$), and a group-theoretic 
{\bf section} $s:G_k\to {\Pi}^{(\et, \ab)}_{\Cal X}$ of the corresponding exact sequence 
$1 \to \Delta ^{\et,\ab}_{\overline {\Cal X}}   \to  {\Pi}^{(\et, \ab)}_{\Cal X} \to    G_k \to 1 $
(cf. diagram (0.4), and the discussion before Theorem A) which is {\bf non-geometric}, i.e., which doesn't arise from a $k$-rational point of $X$.
\endproclaim

\demo {Proof}
Indeed, this follows from Theorem A, Remark 0.2 (i), and the fact that there exists a $p$-adic local field $k$ and a $k$-curve $X$ satisfying the assumptions in Theorem B
such that $X(k)=\emptyset$ (cf. [Sa\"\i di1], Proof of Proposition 3.2.1).
\qed
\enddemo

More generally, we obtain the following.

\proclaim {Theorem C} There {\bf exists} a $p$-adic local field $k$, 
and a smooth, projective, geometrically connected, and {\bf hyperbolic curve} $X$ over
$k$ such that the following holds. The curve $X$ has {\bf bad semi-stable reduction} over $\Cal O_k$,
there exists a {\bf proper} and {\bf flat semi-stable regular} model $\Cal X\to \Spec \Cal O_k$ of $X$ over $\Cal O_k$, and a group-theoretic 
{\bf section} $s:G_k\to {\Pi}^{(\et, \ab)}_{\Cal X}$ of the corresponding exact sequence 
$1 \to \Delta ^{\et,\ab}_{\overline {\Cal X}}   \to  {\Pi}^{(\et, \ab)}_{\Cal X} \to    G_k \to 1 $
(cf. diagram (0.4), and the discussion before Theorem A) which is {\bf non-geometric}, i.e., which doesn't arise from a $k$-rational point of $X$.
\endproclaim

\subhead
\S 1. Group theoretic sections of the \'etale-by-geometrically abelian quotient $\Pi^{(\et,\ab)}$
\endsubhead
We use the notations introduced in $\S0$.
In this section we investigate the group-theoretic splittings of the exact sequence (cf. $\S0$, diagram (0.4))
$$1 \to \Delta ^{\et,\ab}_{\overline {\Cal X}}  \to {\Pi}^{(\et, \ab)}_{\Cal X} \to     G_k \to 1.$$

Let $k$, $X\to \Spec k$, and $\Cal X\to \Spec \Cal O_k$,  be as in the discussion before Theorem A. Thus, 
$$X\to \Spec k$$ is a smooth, proper, and geometrically connected {\it hyperbolic  curve over the $p$-adic local field} $k$, and 
$$\Cal X\to \Spec \Cal O_k$$ 
is a {\it proper}, and {\it flat model} of $X$ over 
the ring of integers $\Cal O_k$ of $k$, such that the {\it $\gcd$ of the multiplicities of the irreducible components of $\Cal X_s$ equals $1$}. 
It follows from the various definitions that
we have a commutative diagram of exact sequences

$$
\CD
@.   1   @.    1 \\
@. @VVV    @VVV \\
@. I_{\overline {\Cal X}}^{\ab}   @= I_{\overline {\Cal X}}^{\ab} \\
@. @VVV    @VVV \\
1 @>>> \Delta ^{\et,\ab}\defeq \Delta ^{\et,\ab}_{\overline {\Cal X}}      @>>>  {\Pi}^{(\et, \ab)}\defeq {\Pi}^{(\et, \ab)}_{\Cal X} @>>>    G_k @>>> 1 \\
@.      @VVV     @VVV   @VVV  \\
1 @>>>  \Delta ^{\et}\defeq  \pi_1(\overline {\Cal X_s},\bar \xi)    @>>>   \Pi ^{(\et)}\defeq \widehat \pi_1(\Cal X_s, \xi)       @>>> G_k  @>>> 1\\
@. @VVV   @VVV \\
@. 1 @. 1 \\
\endCD
\tag 1.1
$$
where we noted
$$\Delta ^{\et}\defeq  \Delta ^{\et}_{\overline {\Cal X}} \defeq\pi_1(\overline {\Cal X_s},\bar \xi), \ \ \ \ \ \Delta ^{\et,\ab}\defeq \Delta ^{\et,\ab}_{\overline {\Cal X}},$$
and
$$\Pi ^{(\et)}\defeq \ \Pi ^{(\et)}_{\Cal X}\defeq \widehat \pi_1(\Cal X_s, \xi),\ \ \ \ \ \  {\Pi}^{(\et, \ab)}\defeq {\Pi}^{(\et, \ab)}_{\Cal X} .$$

The profinite group $\Delta ^{\et}$ is {\bf finitely generated} (as follows from the well-known finite generation of the profinite group 
$\pi_1(\overline X,\overline \eta)$ which projects onto $\Delta ^{\et}$).
Let $\{\Delta ^{\et,i}\}_{i\ge 1}$ be a countable system of {\bf characteristic open} subgroups of $\Delta^{\et}$ such that 
$$\Delta ^{\et,i+1}\subseteq \Delta ^{\et,i},\ \ \ \ \Delta ^{\et,1}\defeq \Delta ^{\et},\ \ \ \ \text {and} \ \ \bigcap _{i\ge 1}\Delta ^{\et,i}=\{1\}.$$
Write $\Delta _i\defeq \Delta ^{\et}/\Delta ^{\et,i}$. Thus, $\Delta _i$ is a {\it finite, characteristic} quotient of $\Delta ^{\et}$, and
we have a {\it pushout} diagram of exact sequences

$$
\CD
1@>>> \Delta^{\et} @>>>   \Pi^{(\et)}  @>>>  G_k @>>> 1\\
@.  @VVV    @VVV    @VVV \\
1@>>>  \Delta _i  @>>>  \Pi _i @>>>  G_k @>>> 1\\
\endCD
\tag 1.2
$$
which defines a ({\it geometrically finite}) quotient $\Pi_i$ of $\Pi^{(\et)}$.

Note that the upper horizontal sequence in diagram (1.2) {\bf splits}. Indeed, the lower horizontal exact sequence in diagram (0.2) splits since the Galois
group $G_F$ of $F$ is pro-free ($F$ being a finite field), hence the middle horizontal sequence in diagram (0.2); which is by definition the sequence
$1\to \Delta ^{\et}\to \Pi^{(\et)}\to G_k\to 1$, splits since the lower right square in that diagram is cartesian. Let   
$$s:G_k\to \Pi^{(\et)}$$
be a {\bf section} of the upper sequence in diagram (1.2), which naturally induces a section
$$s_i:G_k\to \Pi_i$$
of the lower sequence in diagram (1.2), for each $i\ge 1$. Write 
$$\widehat \Pi^i\defeq \widehat \Pi^i[s] \defeq \Delta ^{\et,i}.s (G_k).$$
Thus, $\widehat \Pi^i\subseteq \Pi^{\et}$ is an {\it open} subgroup which contains the image 
$s (G_k)$ of $s$. Write $\Pi^i$ for the {\bf inverse image} of $\widehat \Pi ^i$
in $\pi_1(X,\eta)$. Thus, $\Pi^i\subseteq \pi_1(X,\eta)$ is an {\it open} subgroup which corresponds
to an \'etale cover 
$$X_i\to X_1\defeq X$$
defined over $k$ (since $\Pi^i$ maps onto $G_k$ via the natural projection $\pi_1(X,\eta)\twoheadrightarrow G_k$, by the very 
definition of $\Pi^i$). 
Moreover, it follows from the various definitions that the \'etale cover $X_i\to X$ extends to an \'etale cover
$$\Cal X_i\to \Cal X,$$
defined over $\Cal O_k$.

Note that the \'etale cover $\overline X_i\defeq X_i\times _k\overline k\to \overline X$ is Galois with
Galois group $\Delta _i$, and we have a commutative diagram of \'etale covers

$$
\CD
\overline X_i  @>>>   \overline X \\
@VVV    @VVV \\
X_i  @>>>  X\\
\endCD
$$
where $\overline X_i\to X$ is Galois with Galois group $\Pi_i$, and $\overline X_i\to X_i$ is Galois with Galois group $s_i(G_k)$. 
Moreover, we have a commutative diagram of exact sequences
$$
\CD
@.   1@.    1@.   \\
@. @VVV  @VVV\\
1 @>>>  \Delta ^i= \pi_1(\overline {X_i},\bar \eta)@>>> \Pi^i=\pi_1(X_i,\eta) @>>>  G_k @>>> 1\\
@.   @VVV    @VVV   @|| \\
1 @>>> \pi_1(\overline {X},\bar \eta) @>>>   \pi_1(X,\eta) @>>>  G_k  @>>> 1 \\
\endCD
$$
where $\Delta^i$ is the {\it inverse image} of $\Delta ^{\et,i}$
in $\pi_1(\overline {X},\bar \eta)$, and the equalities $\Delta ^i= \pi_1(\overline {X_i},\bar \eta)$, $\Pi^i=\pi_1(X_i,\eta)$,
are natural identifications; the base points $\eta$ (resp. $\overline \eta$) of $X_i$ (resp. $\overline X_i$) are those induced by 
the base points $\eta$ (resp. $\overline \eta$) of $X$ (resp. $\overline X$). Note that $\Pi^{i+1}\subseteq \Pi^i$ and $\Delta ^{i+1}\subseteq \Delta ^i$ 
as follows from the various definitions.

\proclaim {Lemma 1.1} With the same notations as above, the following holds:
$$I_{\overline {\Cal X}}=\bigcap _{i\ge 1} \Delta ^i,\ \ \ \ \ \ \ \  \text {and}\ \ \ \ \ \ \ \  I_{\Cal X}=\bigcap _{i\ge 1} \Pi ^i.$$
\endproclaim

\demo{Proof}
Follows from the various definitions.
\qed
\enddemo

For each integer $i\ge 1$, write $\Delta^{i, \ab}$ for the maximal {\it abelian} quotient of $\Delta ^i$; which is a {\it characteristic} quotient.  
Consider the natural {\it pushout} diagram
$$
\CD
1 @>>>  \Delta ^i= \pi_1(\overline {X_i},\bar \eta)@>>> \Pi^i=\pi_1(X_i,\eta) @>>>  G_k @>>> 1\\
@.   @VVV    @VVV   @VVV\\
1 @>>>  \Delta ^{i,\ab}@>>> \Pi^{(i,\ab)}@>>>  G_k @>>> 1\\
\endCD
$$
Thus,  $\Pi^{(i,\ab)}$ is the {\bf geometrically abelian} fundamental group of $X_i$.
Write $J_i\defeq J_{X_i}$ for the jacobian variety of $X_i$, and $\overline {J_i}\defeq J_i\times _k\overline k$ for the jacobian variety of $\overline X_i$.
Let $J_i^1$ be the $J_i$-torsor $\Pic^1_{X_i}$. We have $J_i^1(k)\neq \emptyset$ if $X_i$ has {\bf index} $1$. 
In this case we have an identification $J_i\isom J_i^1$. Moreover, $\Delta ^{i,\ab}$ 
is naturally identified with the Tate module $T\overline {J_i}$ of $\overline {J_i}$ as $G_k$-modules.

\proclaim {Lemma 1.2} We use the above notations. The exact sequence
$$1 \to  \Delta ^{i,\ab} \to \Pi^{(i,\ab)} \to  G_k \to 1$$
is a {\bf split} exact sequence of profinite groups {\bf if} the {\bf index} of $X_i$ equals $1$.
\endproclaim

\demo{Proof} Indeed, this exact sequence is naturally identified with the exact sequence $1\to  \pi_1(\overline {X_i},\bar \eta)^{\ab} \to 
\pi_1(X_i,\eta)^{(\ab)} \to  G_k \to 1$, where $\pi_1(\overline {X_i},\bar \eta)^{\ab}$ is the maximal abelian quotient of $\pi_1(\overline {X_i},\bar \eta)$,
and $\pi_1(X_i,\eta)^{(\ab)}$ is the geometrically abelian quotient of $X_i$. Moreover, it is well-known that the exact sequence 
$1\to  \pi_1(\overline {X_i},\bar \eta)^{\ab} \to \pi_1(X_i,\eta)^{(\ab)} \to  G_k \to 1$ splits if the class $[J_i^1]$ of $J_i^1$ in $H^1(G_k,J_i)$
is trivial (cf. [Harari-Szamuely], Theorem 1.2,  for a more general/precise statement).
\qed
\enddemo

Consider the natural {\it pullback} commutative diagram 

$$
\CD
1 @>>> I_{\overline {\Cal X}}^{\ab}  @>>> \Cal H_{\Cal X}\defeq   \Cal H_{\Cal X}[s] @>>>  G_k @>>> 1\\
@.   @VVV    @VVV   @VsVV\\
1 @>>>   I_{\overline {\Cal X}}^{\ab}  @>>> \Pi^{(\et,\ab)}  @>>>  \Pi ^{(\et)} @>>> 1 \\ 
\endCD
\tag 1.3$$ 
 where the right square is cartesian. Thus, (the group extension) $\Cal H_{\Cal X}$ is the {\it pullback} of (the group extension) 
 $\Pi^{(\et,\ab)}$ via the section  $s:G_k\to \Pi^{(\et)}$.

 \proclaim {Lemma 1.3} We have natural identifications $I_{\overline {\Cal X}}^{\ab}  \isom \underset{i\ge 1} \to{\varprojlim}\ \Delta ^{i,\ab}$,
 and
 $\Cal H_{\Cal X}\isom \underset{i\ge 1} \to{\varprojlim}\  \Pi^{(i,\ab)}$.
 \endproclaim
 
 \demo {Proof} Follows from the various definitions. More precisely, recall the \'etale Galois cover $\overline X_i\to \overline X$
 with Galois group $\Delta _i$. For each integer $j\ge 1$, let $\Delta_{i,j}\defeq \Delta ^{i,\ab}/j.\Delta ^{i,\ab}$, which is a
 characteristic quotient of $\Delta ^{i,\ab}$, and $\overline X_{i,j}\to \overline X_i$ the corresponding \'etale abelian cover
 with Galois group $\Delta _{i,j}$.
 Then the \'etale cover $\overline X_{i,j}\to \overline X$ is Galois with Galois group $\widetilde \Delta _{i,j}$, which inserts in
 the following exact sequence $1\to \Delta _{i,j}\to \widetilde \Delta _{i,j}\to \Delta _i\to 1$. The $\{\widetilde \Delta _{i,j}\}_{(i,j)}$,
 where the set of pairs $(i,j)$ is endowed with the product order, form a cofinal system of finite quotients of $\Delta ^{\et,\ab}$.
 From this it follows that the  $\{\Delta _{i,j}\}_{(i,j)}$ form a cofinal system of finite quotients of $I_{\overline {\Cal X}}^{\ab}$.
 Thus, $I_{\overline {\Cal X}}^{\ab}\isom \underset{(i,j)} \to{\varprojlim}\ \Delta _{i,j}\isom 
 \underset{i\ge 1} \to{\varprojlim}\ (\underset{j\ge 1} \to{\varprojlim}\ \Delta _{i,j})\isom \underset{i\ge 1} \to{\varprojlim}\ \Delta ^{i,\ab}$.
 The identification $\Cal H_{\Cal X}\isom \underset{i\ge 1} \to{\varprojlim}\  \Pi^{(i,\ab)}$ follows immediately.
  \qed
 \enddemo

\proclaim {Lemma 1.4} The exact sequence $1\to \Delta ^{\et,\ab}   \to  {\Pi}^{(\et, \ab)} \to    G_k \to 1 $ {\bf splits if}
the exact sequence $1 \to I_{\overline {\Cal X}}^{\ab} \to \Cal H_{\Cal X}  \to G_k \to 1$ {\bf splits}.
\endproclaim

\demo{Proof}
Follows immediately from diagram (1.3). 
\qed
\enddemo

 \proclaim {Proposition 1.5} Suppose that the exact sequence $1 \to  \Delta ^{i,\ab} \to \Pi^{(i,\ab)} \to  G_k \to 1$
 splits (this is the case for example if the curve $X_i$ has index $1$ (cf. Lemma 1.2)), for each integer $i\ge 1$. 
 Then the exact sequence $1 \to I_{\overline {\Cal X}}^{\ab} \to \Cal H_{\Cal X}  \to G_k \to 1$ {\bf splits}.
 \endproclaim
 
 \demo {Proof} We will show the existence of a section $\tilde s:G_k\to \Cal H_{\Cal X}$ of the natural projection 
 $\Cal H_{\Cal X}\twoheadrightarrow G_k$. For each integer $i\ge 1$ let 
 $$\{\Delta^{i,\ab}_j\defeq j.\Delta ^{i,\ab}\}_{j\ge 1},$$ 
 which is a system of open 
 characteristic subgroups of $\Delta ^{i,\ab}$ such that 
$$\Delta ^{i,\ab}_{j'}\subseteq \Delta ^{i,\ab}_j\ \ \  \text {if}\ \ \  j/j',\ \ \ \ \Delta ^{i,\ab}_1\defeq \Delta ^{i,\ab},\ \ \ \ \text {and} \ \ \bigcap _{j\ge 1}\Delta ^{i,\ab}_j=\{1\}.$$
Write $\Delta_{i,j}\defeq \Delta ^{i,\ab}/\Delta ^{i,\ab}_j$, and $\Pi_{i,j}\defeq \Pi ^{(i,\ab)}/\Delta ^{i,\ab}_j$. Thus, we have a natural exact sequence
$$1\to \Delta _{i,j} \to \Pi_{i,j} \to G_k\to 1, $$
and a projective system $\{\Pi_{i,j}\}_{(i,j)}$, where the the set of pairs of integers $(i,j)$ is endowed with the product order.
Note that we have natural identifications $\Delta ^{(i,\ab)} \isom \underset{j\ge 1} \to{\varprojlim}\  \Delta_{i,j}$, and
$\Pi ^{(i,\ab)} \isom \underset{j\ge 1} \to{\varprojlim}\  \Pi_{i,j}$.
Moreover, we have a natural identification $\Cal H_{\Cal X}\isom \underset{(i,j)} \to{\varprojlim}\  \Pi_{i,j}$ (cf. Lemma 1.3, and the above discussion). 
In particular, the set $\Sec (G_k,\Cal H_{\Cal X})$ of group-theoretic sections of the natural projection $\Cal H_{\Cal X}\twoheadrightarrow G_k$ is naturally identified with
the projective limit $\underset{(i,j)} \to{\varprojlim}\  \Sec(G_k,\Pi_{i,j})$ of the sets  $\Sec(G_k,\Pi_{i,j})$
of group-theoretic sections of the natural projection 
$\Pi_{i,j}\twoheadrightarrow G_k$. For each pair of integers $(i,j)$ the set $\Sec(G_k,\Pi_{i,j})$ is non-empty, since $\Pi_{i,j}$ is a quotient
of $\Pi^{(i,\ab)}$ (cf. the assumption that (the group extension) $\Pi^{(i,\ab)}$ splits). Moreover, the set $\Sec(G_k,\Pi_{i,j})$ is, up to conjugation by the elements of $\Delta_{(i,j)}$, a torsor under the group
$H^1(G_k,\Delta _{i,j})$ which is finite since $k$ is a $p$-adic local field (cf. [Neukirch-Schmidt-Winberg], (7.1.8) Theorem (iii)). 
Thus, $\Sec(G_k,\Pi_{i,j})$ is a nonempty finite set.
Hence the set $\Sec (G_k,\Cal H_{\Cal X})$ is nonempty being the projective limit of nonempty finite sets.
 This finishes the proof of Proposition 1.5.
 \qed 
 \enddemo
 
 In fact, we proved the following more precise statement.
 
\proclaim {Proposition 1.6} Let $s:G_k\to \Pi_{\Cal X}^{(\et)}$ be a section 
of the natural projection $\Pi_{\Cal X}^{(\et)}\twoheadrightarrow G_k$. Suppose that for each
open subgroup $H\subseteq \Pi_{\Cal X}^{(\et)}$ such that $s(G_k)\subset H$, and the
corresponding \'etale cover $\Cal X_H\to \Cal X$, it holds that the index of 
the $k$-curve $\Cal X_H\times _{\Cal O_k}k$ equals 1.
Then there exists a section
 $\tilde s:G_k\to \Pi_{\Cal X}^{(\et,\ab)}$  
 of the natural projection $\Pi_{\Cal X}^{(\et,\ab)}\twoheadrightarrow G_k$ which lifts the section $s$, i.e., such that we have 
 a commutative diagram
 $$
 \CD
 G_k  @>{\tilde s}>> \Pi_{\Cal X}^{(\et,\ab)}\\
 @V{\id}VV         @VVV\\
 G_k  @>{s}>>  \Pi_{\Cal X}^{(\et)}\\
 \endCD
 $$
 where the right vertical map is the map in diagram (1.1).
 \endproclaim
 
 \demo {Proof} Follows from the above discussion, Proposition 1.4, Proposition 1.5, the commutative diagram (1.3), and Lemma 1.2.
 \qed
 \enddemo
 
 \definition {Remark 1.7} Proposition 1.6 remains valid if instead of assuming the curves $\Cal X_H\times _{\Cal O_k}k$ to have index $1$ one assumes that
 if $J_H$ denotes the jacobian of   $\Cal X_H\times _{\Cal O_k}k$ and $J_H^1$ is the corresponding torsor $\Pic^1$, then the class $[J_H^1]$ of $J_H^1$ in
 $H^1(G_k,J_H)$ lies in the maximal divisible subgroup of $H^1(G_k,J_H)$, as Lemma 1.2 is still valid under a similar assertion 
 (cf. proof of Lemma 1.2 and the reference therein).
 \enddefinition

\subhead
\S 2. Proof of Theorem A and Theorem C
\endsubhead
The rest of this paper is devoted to proving Theorem A and Theorem C. We use the notations introduced in $\S0$ and $\S1$.

\demo {Proof of Theorem A}
Assume that $\Cal X$ satisfies the condition $(\star)$.
Thus, $\Cal X$ is {\it regular}, the special fibre $\Cal X_s\defeq \Cal X\times_{\Spec \Cal O_k} \Spec F$ has an {\it irreducible component} $X_{i_0}$
which is {\it (geometrically) reduced}, {\it geometrically irreducible}, {\it geometrically unibranch}, and {\it its normalisation $X_{i_0}^{\nor}$ is hyperbolic}.  
In order to prove Theorem A it suffices to construct a group-theoretic section 
$s:G_k\to \Pi^{(\et)}$ of the natural projection $\Pi^{(\et)}\twoheadrightarrow G_k$, with corresponding open subgroups $\widehat \Pi^i[s]\subset \Pi^{(\et)}$,
such that in the corresponding \'etale covers $X_i\to X$ the index of $X_i$ equals $1$, for each integer $i\ge 1$ (cf. Proposition 1.5 and Lemma 1.4).

Let 
$$D_{X_{i_0}}\subset \Pi^{(\et)}\defeq \widehat \pi_1(\Cal X_s, \xi)$$ 
be a {\it decomposition group} associated to the irreducible component $X_{i_0}$ of $\Cal X_s$. More precisely, $D_{X_{i_0}}$ is the decomposition 
group of a connected component $\widetilde X_{i_0}$ of the fibre of $X_{i_0}$  in (the special fibre of) the universal pro-\'etale cover of $\Cal X$, 
corresponding to the quotient $\Pi^{(\et)}$ of $\pi_1(X,\eta)$,
and $D_{X_{i_0}}$ is only defined up to conjugation.  Note that $\widetilde X_{i_0}$ is an {\it irreducible pro-curve} (i.e., is the projective limit of 
irreducible curves) since the morphism $\widetilde X_{i_0}\to X_{i_0}$ is pro-\'etale, $X_{i_0}$ is geometrically unibranch, and $\widetilde X_{i_0}$ is connected.

We have a commutative diagram
$$
\CD
1 @>>>   D_{\overline X_{i_0}}  @>>>   D_{X_{i_0}}  @>>>   G_k  @>>>  1\\
@.  @VVV          @VVV          @VVV \\
1  @>>> \pi_1(\overline X_{i_0},\bar \xi)  @>>>  \pi_1 (X_{i_0}, \xi)   @>>>  G_F  @>>> \\
\endCD
\tag 2.1$$
where $D_{\overline X_{i_0}}$ is defined so that the upper sequence is exact, 
$\pi_1 (X_{i_0}, \xi)$ is the {\it arithmetic fundamental group} of $X_{i_0}$, the right vertical map is the natural projection, and the right square is {\it cartesian}. 
In particular, the left vertical map in the above diagram is an isomorphism. 

Let 
$$\bar s_{i_0}:G_F\to \pi_1 (X_{i_0}, \xi)$$ 
be a group-theoretic section of the lower sequence in diagram (2.1) (this sequence splits since  $G_F$ is pro-free), 
which induces naturally (cf. above cartesian square in diagram (2.1)) a section
$$s_{i_0}:G_k\to D_{X_{i_0}}$$ 
of the natural projection 
$D_{X_{i_0}}\twoheadrightarrow G_k$, hence induces also a {\bf section}  
$$s\defeq s_{i_0} :G_k\to D_{X_{i_0}}\subset \Pi^{(\et)}$$
of the upper sequence in diagram (1.2) (cf. above discussion). 

Recall the above \'etale cover
$\Cal X_i\to \Cal X$ (cf. discussion after diagram (1.2)), which corresponds to the open subgroup $\widehat \Pi^i[s]=\Delta ^{i,\et}.s(G_k)$.
Let $X_{i,i_0}$ be the image of $\widetilde X_{i_0}$ in $\Cal X_i$, and 
$X_{i,i_o}\to X_{i_0}$ the corresponding \'etale cover. Thus, the \'etale cover $X_{i,i_o}\to X_{i_0}$ corresponds to  
the open subgroup $\widehat \Pi^i[s]\cap D_{X_{i_0}}\subseteq D_{X_{i_0}}$ of $D_{X_{i_0}}$.
Then $X_{i,i_0}$ is an {\it $F$-curve} and the cover $X_{i,i_0}\to X_{i_0}$ is a morphism of $F$-curves, since the open subgroup 
$\widehat \Pi^i[s]\cap D_{X_{i_0}}\subseteq D_{X_{i_0}}$ contains $s_{i_0}(G_k)$ (hence projects onto $G_k$ via the projection
$D_{X_{i_0}}\twoheadrightarrow G_k$), $X_{i,i_0}$ is
{\it reduced} (since it is an \'etale cover of $X_{i_0}$ which is reduced), hence is also geometrically reduced since 
$F$ is perfect, and $X_{i,i_0}$ is {\it geometrically irreducible} (since $\widetilde X_{i_0}$ is irreducible (cf. above discussion)). More precisely, 
$X_{i,i_0}$ is geometrically unibranch; being an \'etale cover of $X_{i_0}$ which is geometrically  unibranch, and is geometrically connected. Moreover, 
its normalisation $X_{i,i_0}^{\nor}$ is {\it hyperbolic}, since $X_{i,i_0}^{\nor}$ dominates $X_{i_0}^{\nor}$.
In particular, the {\bf index} of $X_i\defeq \Cal {X}_i\times _{\Spec \Cal O_k}\Spec k$ equals $1$ (cf. Remark 0.2 (ii)).

This finishes the proof of Theorem A.
\qed
\enddemo

\demo{Proof of Theorem C} Let $F$ be a finite field, and $C\to \Spec F$ a singular {\bf stable} curve with arithmetic genus
$g(C)>1$ such that each double point of $C$ is $F$-{\bf rational} and lies on two distinct irreducible components of $C$. For example let $X_0\to \Spec F$
be a proper, smooth, and geometrically connected hyperbolic curve such that $X_0(F)=\{x_i\}_{i\in I}\neq \emptyset$. For each $i\in I$, let $E_i$ be an elliptic curve 
over $F$, and identify the origin of $E_i$ with the rational point $x_i$ of $X_0$ into an ordinary $F$-rational double point $x_i$. 
We obtain a singular stable curve $C=X_0+\sum _{i\in I}E_i$, whose configuration is tree like, where $X_0$ intersects $E_i$ at the rational double point $x_i$,
and $E_i\cap E_j=\emptyset$, for $i \neq j$. Let $k$ be a $p$-adic local field with residue field $F$ and ring of integers $\Cal O_k$. Using formal patching techniques
one can construct a proper and {\bf stable} relative curve $\Cal X\to \Spec \Cal O_k$ such that $\Cal X$ is {\bf regular}, and $\Cal X\times _{\Cal O_k}F=C$
(compare with [Sa\"\i di2], Proposition 3.6). Let $X\defeq \Cal X\times _{\Cal O_k}k$. The morphism $x_i:\Spec F\to C$ gives rise naturally to a (conjugacy class of)
section $\bar s_{x_i}:G_F\to \pi_1(C)$ of the natural projection 
$\pi_1(C)\twoheadrightarrow G_F$, which naturally induces a section $s:G_k\to \Pi^{(\et)}\defeq \Pi^{(\et)}_{\Cal X}$ 
of the natural projection $\Pi^{(\et)}\twoheadrightarrow G_k$ (cf. diagram 0.2).
Let $H\subseteq \Pi^{(\et)}$ be an open subgroup such that $s(G_k)\subset H$, and $\Cal X_H\to \Cal X$ the corresponding 
\'etale cover. Thus, $\Cal X_H\to \Spec \Cal O_k$ is a relative stable curve, and $\Cal X_H$ is regular.  
Let $X_H\defeq \Cal X_H\times _{\Cal O_k}k$, then the {\bf index} of $X_H$ equals $1$. 
Indeed, if $C'$ is an irreducible component of $C$ which contain the (double) point $x_i$, and $C_H'$ is an irreducible component 
of $\Cal X_H$ above $C'$ containing an $F$-rational (double) point above $x_i$ (cf. the construction of the section $s$), 
then $C_H'$ is geometrically integral (cf. Remark 0.2 (ii)).

The above section $s$ satisfies the assumption of Proposition 1.6.
Thus, there exists a section $\tilde s:G_k\to \Pi^{(\et,\ab)}$ of the natural projection $\Pi^{(\et,\ab)}\twoheadrightarrow G_k$ which lifts the section $s$ (cf. loc. cit.).
The section $\tilde s$ is {\bf not geometric}. Indeed, if $\tilde s$ is geometric and arises from a rational point $x\in X(k)$, then $s$ is also geometric and arises from the 
rational point $x$. Then, by construction of the section $s$, the point
$x$ would specialise in the double point $x_i$ as is easily verified. But this is not possible since $\Cal X$ is regular and $x$ should specialise in a smooth point of $C$
(cf. [Liu], Corollary 9.1.32). 

This finishes the proof of Theorem C. 
\qed
\enddemo

 \definition {Acknowledgment} I would like very much to thank Akio Tamagawa for several discussions we had on the topic of this paper. 
 I also thank the referee for his/her valuable comments.
 \enddefinition

$$\text{References.}$$

\noindent
[Grothendieck] Grothendieck, A., Rev\^etements \'etales et groupe fondamental, Lecture 
Notes in Math. 224, Springer, Heidelberg, 1971.

\noindent
[Harari-Szamuely] Harari, D., Szamuley, T., Galois sections for abelianized fundamental groups (With an appendix by E. V. Flynn), 
Math. Ann. 344 (2009), no. 4, 779--800.

\noindent
[Hoshi] Hoshi, Y., Existence of nongeometric pro-$p$ Galois sections of hyperbolic curves, Publ. Res. Inst. Math. Sci. 46 (2010), no. 4, 829-848.

\noindent
[Liu] Liu, Q., Algebraic geometry and arithmetic curves, Oxford graduate texts in mathematics 6. Oxford University Press, 2002.

\noindent
[Neukirch-Schmidt-Winberg] Neukirch, J., Schmidt, A., Winberg, K., Cohomology of Number Fields, Grundlehren der mathematischen Wissenschaften, 323, Springer, 2000.

\noindent
[Raynaud] Raynaud, M., Sp\'ecialisation du fonctuer de Picard, Publications math\'ematiques de l'I.H.\'E.S., tome 38 (1970), p. 27-76.

\noindent
[Sa\"\i di] Sa\"\i di, M., On the $p$-adic section conjecture, Journal of Pure and Applied Algebra, Volume 217, Issue 3 (2013), 583-584.

\noindent
[Sa\"\i di1] Sa\"\i di, M., The cuspidalisation of sections of arithmetic fundamental groups, Advances in Mathematics 230
(2012) 1931-1954.

\noindent
[Sa\"\i di2] Sa\"\i di, M., Rev\^etements \'etale ab\'eliens, courants sur les graphes et 
r\'eduction semi-stable des courbes, Manuscripta Math. 89 (1996), no.2, 245-265.

\bigskip

\noindent
Mohamed Sa\"\i di

\noindent
College of Engineering, Mathematics, and Physical Sciences

\noindent
University of Exeter

\noindent
Harrison Building

\noindent
North Park Road

\noindent
EXETER EX4 4QF

\noindent
United Kingdom

\noindent
M.Saidi\@exeter.ac.uk

\end
\enddocument